\documentclass{amsart}
\usepackage{amsthm,amssymb}

\newtheorem{proposition}{Proposition}[section]
\newtheorem{theorem}[proposition]{Theorem}
\newtheorem{corollary}[proposition]{Corollary}
\newtheorem{lemma}[proposition]{Lemma}
\newtheorem{conjecture}[proposition]{Conjecture}

\theoremstyle{definition}
\newtheorem{remark}[proposition]{Remark}
\newtheorem*{definition}{Definition}

\newtheorem{example}[proposition]{Example}

\newcommand{\itref}[1]{(\ref{#1})}

\newcommand{\re}{||}
\newcommand{\bdmap}{\partial}
\newcommand{\calC}{\mathcal C}

\newcommand{\calS}{\mathcal S}
\newcommand{\calN}{\mathcal N}
\newcommand{\IN}{IN}
\newcommand{\s}{\mathbf{s}}
\newcommand{\ostar}{\circ}
\newcommand{\betat}{\tilde{\beta}}
\newcommand{\Ht}{\tilde{H}}
\newcommand{\abs}[1]{\lvert#1\rvert}
\newcommand{\setm}{\backslash}
\newcommand{\skel}[1]{^{(#1)}}
\newcommand{\reals}{\mathbb R}
\newcommand{\disun}{\mathbin{\dot{\cup}}} 
\newcommand{\condns}[2]{\substack{#1 \\ #2}}
\newcommand{\blambda}{\boldsymbol{\lambda}}
\newcommand{\bmu}{\boldsymbol{\mu}}
\DeclareMathOperator{\rk}{rank}
\DeclareMathOperator{\ci}{ci}
\DeclareMathOperator{\im}{im}
\DeclareMathOperator{\st}{st}

\newcommand{\ie}{{\em i.e.}}
\newcommand{\eg}{{\em e.g.}}

\begin{document}

\title{A Relative Laplacian spectral recursion}

\author{Art M. Duval}
\email{artduval@math.utep.edu}
\address{Department of Mathematical Sciences\\
         University of Texas at El Paso\\
         El Paso, TX 79968-0514}



\begin{abstract}

The Laplacian spectral recursion, satisfied by matroid complexes and
shifted complexes, expresses the eigenvalues of the combinatorial
Laplacian of a simplicial complex in terms of its deletion and
contraction with respect to vertex $e$, and the relative simplicial
pair of the deletion modulo the contraction.  We generalize this
recursion to relative simplicial pairs, which we interpret as
intervals in the Boolean algebra.  The deletion modulo contraction
term is replaced by the result of removing from the interval all pairs
of faces in the interval that differ only by vertex $e$.

We show that shifted pairs and some matroid pairs satisfy this recursion.
We also show that the class of intervals satisfying this recursion is
closed under a wide variety of operations, including duality and
taking skeleta.
\end{abstract}

\maketitle

\section{Introduction}\label{se:intro}

There are two good reasons to extend the Laplacian spectral recursion
from simplicial complexes to relative simplicial pairs.

The spectral recursion for simplicial complexes expresses the
eigenvalues of the combinatorial Laplacian $\bdmap \bdmap^* + \bdmap^*
\bdmap$ of a simplicial complex $\Delta$ in terms of the eigenvalues
of its deletion $\Delta-e$, contraction $\Delta/e$, and an ``error
term'' $(\Delta-e,\Delta/e)$.  This recursion does not hold for all
simplicial complexes, but does hold for independence complexes of
matroids and shifted simplicial complexes \cite{Duval:recursion}.  In
each case, the deletion and contraction are again matroids or shifted
complexes, respectively, but the error term is only a relative
simplicial pair of the appropriate kind of complexes.  Being able to
apply the recursion to relative simplicial pairs, such as the error
term, would make the spectral recursion truly recursive.

A more compelling reason 
comes from duality, the idea that a Boolean algebra looks the same
upside-down as it does right-side-up.  
Many operations preserve the property of satisfying the spectral
recursion \cite{Duval:recursion}, but the {\em dual} $\Delta^*$ (see
equation \eqref{eq:dual}) of simplicial complex $\Delta$, which is an
order filter instead of a simplicial complex, satisfies only a
slightly modified version of the spectral recursion when $\Delta$
satisfies the spectral recursion \cite[Theorem 6.3]{Duval:recursion}.
Relative simplicial pairs include both simplicial complexes and order
filters as special cases, and so suggest a way to unify the two
versions of the spectral recursion.

Furthermore, the Laplacian itself is self-dual (Section
\ref{se:laplacians}), and so we will state and prove most of our
results in self-dual form.  The first step is to think of relative
simplicial pairs as {\em intervals} in the Boolean algebra of subsets
of the set of vertices, since the dual of an interval is again an
interval, in a very natural way.  To further emphasize this symmetry,
we represent these intervals by vertically symmetric capital Greek
letters, such as $\Phi$ and $\Theta$.
When we extend the spectral recursion from simplicial complexes to
intervals, the ideas of deletion and contraction generalize easily and
naturally.  But, even with duality as a guide, it is not as clear
what should replace $(\Delta-e,\Delta/e)$ as the error term.

The answer turns out to be to remove from $\Phi$ all the pairs $\{F,F
\disun e\}$ in $\Phi$.  This simple operation, which we will call the
{\em reduction} of interval $\Phi$ with respect to $e$, and denote by
$\Phi \re e$, has a few remarkable (but easy to prove) properties that
will allow us to show that it is the correct error term.  To start, it
is clear that this operation is self-dual, which goes nicely with
deletion and contraction being more or less duals of one another.
Somewhat more surprising is that $\Phi \re e$ is still an interval,
albeit in two separate components (Lemma \ref{th:par.interval} and
Proposition \ref{th:par.direct.sum}).  Finally, it is necessary for
the error term to have the same homology as $\Phi$ itself (see Lemma
\ref{th:recursion.special.values}), and $\Phi \re e$ satisfies this as
well (equation \eqref{eq:beta.m0}).  Perhaps reduction deserves
further investigation, beyond Laplacians, since it is easy to compute,
preserves homology, and produces a smaller interval.  (Reduction is a
special case of collapsing induced by a discrete Morse function coming
from an acyclic, or Morse, matching, $F \leftrightarrow F \disun e$,
for all possible $F$; see \cite{Chari,Forman}.)

Of course, the most important evidence that reduction is the right
answer is that the spectral recursion for intervals, with $\Phi \re e$
as the error term (equation \eqref{eq:big}), holds for a variety of
intervals.  We are able to prove (Theorem \ref{th:shifted.big}) that
it does hold for shifted intervals, that is, relative simplicial pairs
of complexes, each of which is shifted on the same ordered vertex set.
The analogue for matroids would be relative simplicial pairs of
matroids connected by a strong map, and here our success is more
limited.  Although experimental evidence supports the conjecture that
the spectral recursion holds for all such pairs (Conjecture
\ref{th:all.matroids}), we are only able to prove it in the case where
the difference in ranks between the matroids is 1 (Theorem
\ref{th:matroid.big}).  This does at least provide strong evidence
that $\Phi \re e$ is the correct error term.  Further evidence is that
the property of satisfying the spectral recursion is closed under many
operations on intervals (Section \ref{se:laplacians}), including
duality (Proposition \ref{th:recursion.dual}).

We formally define intervals and their operations, including
reduction, in Section \ref{se:intervals}.  We review Laplacians and
introduce the spectral recursion for intervals in Section
\ref{se:laplacians}.  Our main results, that skeleta preserve the
property of satisfying the spectral recursion (Theorem
\ref{th:recursion.skeleta.hard}), and that shifted intervals and
certain matroid pairs satisfy the spectral recursion (Theorems
\ref{th:shifted.big} and \ref{th:matroid.big}), are the foci of
Sections \ref{se:skeleta}, \ref{se:shifted}, and \ref{se:matroid},
respectively.

\section{Intervals}\label{se:intervals}

In this section, we formally define intervals, and extend many
simplicial complex operations to intervals.  We also introduce the
reduction operation ($\Phi \re e$), and establish some of its
properties.

\begin{definition} 
Let $2^E$ denote the Boolean algebra of subsets of finite set $E$.  We
will say $\Phi \subseteq 2^E$ is an {\em interval} if $F \subseteq G
\subseteq H$ and $F,H \in \Phi$ together imply $G \in \Phi$.  We will
call the set $E$ the {\em ground set} of $\Phi$, individual members of
$E$ the {\em vertices} of $\Phi$, and members of $\Phi$ the {\em
faces} of $\Phi$.
Note that $v$ may be a vertex of interval $\Phi$ without being in {\em
any} face of $\Phi$.  In this case we call $v$ a {\em loop} of $\Phi$.
(This is in analogy to a loop of a matroid.)
\end{definition}

An ``interval'' could be similarly defined on any partially ordered
set, not just $(2^E, \subseteq)$.  Indeed, later on (Section
\ref{se:shifted}), we will consider ``intervals'' on $2^E$ with
respect to a different partial order.  But what makes Laplacians work
so well on intervals of $(2^E, \subseteq)$ is that $(2^E,\subseteq)$
forms a chain complex (Lemma \ref{th:betat.par}, and the preceding
disucssion).  Hereinafter, the word ``interval'' will only refer to
intervals on $(2^E,\subseteq)$.

An important special case of an interval is a {\em simplicial
complex}.  As usual, $\Delta \subseteq 2^E$ is a simplicial complex if
$G \subseteq H$ and $H \in \Delta$ together imply $G \in \Delta$.  It
is obvious that simplicial complexes may be defined as intervals
containing the empty face $\emptyset$.  Of course, our motivation runs
in the oppposite direction; intervals are usually presented as pairs
of simplicial complexes.  If $\Delta' \subseteq \Delta$ are a pair of
simplicial complexes on the same vertex set, then the {\em relative
simplicial pair} $(\Delta,\Delta')$ is simply the set difference
$\Delta - \Delta'$.  We now formally check that intervals and
relative simplicial pairs represent the same objects.

\begin{lemma}\label{th:interval.rel.pair}
Let $\Phi \subseteq 2^E$.  Then $\Phi$ is an interval iff
$\Phi=(\Delta,\Delta')$ for some simplicial complexes $\Delta,\Delta'$.
\end{lemma}
\begin{proof}
To prove the backwards implication, assume $F \subseteq G \subseteq H$
and $F,H \in (\Delta,\Delta')$, so $F,H \in \Delta$, but $F,H \not\in
\Delta'$.  From $G \subseteq H \in \Delta$, we conclude $G \in
\Delta$, but from $F \subseteq G$ and $F \not\in \Delta'$ we conclude
$G \not\in \Delta'$.  Thus $G \in (\Delta,\Delta')$ as desired.

To prove the forwards implication, let $\Delta=\{G \subseteq E\colon G
\subseteq H\ {\rm for\ some}\ H \in \Phi\}$, and let $\Delta' = \Delta
- \Phi$.  Now, if $F \subseteq G \in \Delta$, then $F \subseteq G
\subseteq H$ for some $H \in \Phi$, so $F \in \Delta$, and thus
$\Delta$ is a simplicial complex.  If $F \subseteq G \in \Delta'$,
then $F \in \Delta$, since $F \subseteq G \in \Delta$; but $F \not\in
\Phi$, since $G \not\in \Phi$ and $F \subseteq G \subseteq H$ for some
$H \in \Phi$.  Thus $F \in \Delta'$, and so $\Delta'$ is a simplicial
complex.
\end{proof}

\begin{example}\label{ex:intro}
Let $\Phi$ be the interval on vertex set $\{1,\ldots,6\}$ whose faces
are $\{12456,1245,1246, 1356, 124, 135, 136\}$.  (Here, we are
omitting brackets and commas on individual faces, for clarity.)  It is
easy to check that $\Phi$ is an interval (see also Example
\ref{ex:direct.sum}).  The formulas in Lemma
\ref{th:interval.rel.pair} set $\Delta$ to be the simplicial complex
with facets (maximal faces) $\{12456,1356\}$, and $\Delta'$ the
simplicial complex with facets $\{1256,1456,2456,356,13\}$.  But we
could add the face 34 to both $\Delta$ and $\Delta'$, and they would
still be simplicial complexes such that $\Phi=(\Delta,\Delta')$.
\end{example}

Although Lemma \ref{th:interval.rel.pair} shows that intervals are the
same as relative simplicial pairs, we will strive to put all of our results in
the language of intervals rather than relative simplicial pairs.  One reason is
the potential difficulty in describing properties of the interval in
terms of the pair of simplicial complexes which are not necessarily
unique, as demonstrated in Example \ref{ex:intro}.  Another, as
alluded to in the Introduction, is to better take advantage of
duality.  The {\em dual} of interval $\Phi$ on ground set $E$ is
\begin{equation}\label{eq:dual}
\Phi^* = \{E - F \colon F \in \Phi\}.
\end{equation}
It is easy to see that the dual of an interval is again an interval,
and that $\Phi^{**}=\Phi$.

It is also easy to see the intersection of two intervals is again an
interval, but we have to be more careful with union, even with
disjoint union.  If $\Phi$ and $\Theta$ are disjoint intervals with
faces $F \in \Phi$ and $G \in \Theta$ such that $F \subseteq G$, then
$\Phi \disun \Theta$, the disjoint union of $\Phi$ and $\Theta$, might
not be an interval.  We thus define two intervals $\Phi$ and $\Theta$
to be {\em totally unrelated} if $F \not\subseteq G$ and $G
\not\subseteq F$ whenever $F \in \Phi$ and $G \in \Theta$, and, in
this case, define the {\em direct sum} of $\Phi$ and $\Theta$ to be
$\Phi \oplus \Theta = \Phi \disun \Theta$.  It is easy to check that
the direct sum of two intervals is again an interval.

\begin{example}\label{ex:direct.sum}
It is easy to see that the interval $\Phi$ of Example \ref{ex:intro}
is a direct sum $\{12456,1245,1246,124\} \oplus \{1356,135,136\}$.
The components of the direct sum are indeed totally unrelated, even
though they share many vertices.
\end{example}

The {\em join} of two intervals $\Phi$ and $\Theta$ on disjoint vertex sets is 
$$
\Phi * \Theta = \{F \disun G\colon F \in \Phi, G \in \Theta\}.
$$
When $\Phi$ and $\Theta$ are simplicial complexes, this matches the
usual definition of join.  It is easy to see that the join of two
intervals is again an interval.
Some special cases of the join deserve particular attention.
If $\Phi$ is an interval and $R$ is a set disjoint from the vertices
of $\Phi$, then define
$$
R \ostar \Phi = \{R\} * \Phi = \{R \disun F\colon F \in \Phi\},
$$
the join of $\Phi$ with the interval whose only face is $R$.
If $v$ a vertex not in $\Phi$, then the {\em cone} of $\Phi$ is 
$$
v * \Phi = \{v,\emptyset\} * \Phi,
$$
the join of $\Phi$ with the interval whose two faces are $v$ and the
empty face.  The {\em open star} of $\Phi$ is $v \ostar \Phi$.
Note that 
$$
v * \Phi = \Phi \disun (v \ostar \Phi).
$$ 

Deletion and contraction are well-known concepts from matroid theory,
and were easily extended to simplicial complexes in
\cite{Duval:recursion}.  Now we further extend to intervals.
If $\Phi$ is an interval and $e$ is a vertex of $\Phi$, then the {\em
deletion} and {\em contraction} of $\Phi$ by $e$ are, respectively,
\begin{align*}
\Phi-e &= \{F \in \Phi\colon e \not\in F\};\\
\Phi/e &= \{F-e\colon F \in \Phi, e \in F\}.
\end{align*}
As opposed to the simplicial complex case, $\Phi/e$ is not necessarily
a subset of $\Phi-e$.  As with simplicial complexes, neither $\Phi/e$ nor
$\Phi-e$ contains $e$ in any of their faces, though we stil consider $e$ to
a vertex, albeit a loop, in each case.  It is also easy to check that
$\Phi-e$ and $\Phi/e$ are intervals when $\Phi$ is an interval.
Note that
\begin{align*}
(\Phi-e)^* &= \{E-F\colon F \in \Phi, e\not\in F\} = \{E-F\colon F \in \Phi, e \in E - F\}\\
           &= e \ostar (\Phi^*/e)\\
\intertext{and, similarly,}
(\Phi/e)^* &= \{E-(F-e)\colon F \in \Phi, e \in F\} = \{(E-F) \disun e\colon F \in \Phi, e \not\in E - F\}\\
           &= e \ostar (\Phi^*-e).
\end{align*}

We are now ready to define reduction, which will be a focal point for
most of the rest of our work.
\begin{definition}
If $\Phi$ is an interval and $e$ is a vertex of $\Phi$, then the {\em
star} of $e$ in $\Phi$ is
$$
\st_{\Phi}e = \bigcup_{F, F \disun e \in \Phi} \{F, F \disun e\}
            = e * ((\Phi-e) \cap (\Phi/e)),
$$
and the {\em reduction} of $\Phi$ by $e$ is
$$
\Phi\re e = \Phi - \st_{\Phi}e.
$$
\end{definition}
When $\Phi$ is a simplicial complex, $\st_{\Phi}e$ matches the usual
definition.  It is easy to check that $\st_{\Phi}e$ is an
interval when $\Phi$ is an interval, but $\Phi \re e$ takes a little more work.

\begin{lemma}\label{th:par.interval}
If $\Phi$ is an interval with vertex $e$, then $\Phi \re e$ is again
an interval.
\end{lemma}
\begin{proof}
Assume otherwise, so $F \subseteq G \subseteq H$, and $F,H \in \Phi
\re e$, but $G \not\in \Phi \re e$.  Thus $F,H \in \Phi$, and, since
$\Phi$ is an interval, $G \in \Phi$.

If $e \not\in G$, then $e\not\in F$, and then $F \subseteq F \disun e
\subseteq G \disun e$.  But also $G \not\in \Phi \re e$ implies
$G \disun e \in \Phi$.  Then, since $\Phi$ is an interval, $F \disun e
\in \Phi$, which contradicts $F \in \Phi \re e$.

Similarly, if instead $e\in G$, then $e \in H$, and then $G-e
\subseteq H-e \subseteq H$.  But also $G \not\in \Phi \re e$
implies $G-e \in \Phi$.  Then since $\Phi$ is an interval, $H-e \in
\Phi$, which contradicts $H \in \Phi \re e$.
\end{proof}

\begin{proposition}\label{th:par.direct.sum}
If $\Phi$ is an interval with vertex $e$, then $\Phi \re e$ is the
direct sum of 
$\{F \in \Phi \re e\colon e \not\in F\} = \{F \in \Phi\colon e \not\in F, F \disun e \not\in \Phi\}$ and 
$\{G \in \Phi \re e\colon e     \in G\} = \{G \in \Phi\colon e     \in G, G    -   e \not\in \Phi\}$.
\end{proposition}
\begin{proof}
To show $\Phi \re e$ is the desired direct sum, let $F, G \in \Phi \re
e$ such that $e \not\in F$, and $e \in G$; we must show $F$ and $G$
are unrelated.  Since $e \in G \setm F$, we know $G \not\subseteq F$,
so assume $F \subseteq G$.  Then $F \subseteq F \disun e \subseteq G$.
Since $F, G \in \Phi$, then also $F \disun e \in \Phi$, which
contradicts $F \in \Phi \re e$.
\end{proof}

\begin{example}\label{ex:par}
Let $\Theta$ be the interval of all faces $F \subseteq \{1,\ldots,6\}$
such that $F$ is a subset of $12356$ or $12456$, but also a superset
of $12$, $135$, or $136$.  It is not hard to check that $\Theta \re 3$
is the interval $\Phi$ of Examples \ref{ex:intro} and
\ref{ex:direct.sum}.  The direct sum decomposition of $\Phi=\Theta \re
3$ given in Example \ref{ex:direct.sum} is the one guaranteed by Proposition
\ref{th:par.direct.sum}.
\end{example}

In the special case where $\Phi$ is a simplicial complex, $\{G \in
\Phi \re e\colon e \in G\}$ is empty and $\Phi \re e = (\Phi -e,
\Phi/e)$.
It is easy to check that $(\st_{\Phi} e)^* = \st_{(\Phi^*)} e$, and so
$(\Phi \re e)^* = \Phi^* \re e$.

We review our notation for boundary maps and homology groups of
simplicial complexes (as in \eg, \cite[Chapter 1]{Munkres}).  As
usual, let $\Phi_i$ denote the set of $i$-dimensional faces of $\Phi$,
and let $C_i = C_i(\Phi;\reals) :=
C_i(\Delta;\reals)/C_i(\Delta';\reals)$ denote the $i$-dimensional
oriented $\reals$-chains of $\Phi=(\Delta,\Delta')$, \ie, the formal
$\reals$-linear sums of oriented $i$-dimensional faces $[F]$ such that
$F \in \Phi_i$.  Let $\bdmap_{\Phi;i}=\bdmap_i\colon C_i \rightarrow
C_{i-1}$ denote the usual (signed) {\em boundary operator}.  Via the
natural bases $\Phi_i$ and $\Phi_{i-1}$ for $C_i(\Phi;\reals)$ and
$C_{i-1}(\Phi;\reals)$, respectively, the boundary operator $\bdmap_i$
has an adjoint map called the {\em coboundary operator},
$\bdmap^*_i\colon C_{i-1}(\Phi; \reals) \rightarrow C_{i}(\Phi;
\reals)$; \ie, the matrices representing $\bdmap$ and $\bdmap^*$ in
the natural bases are transposes of one another.

As long as $\Phi$ is an interval, $C_{\bullet}(\Phi;\reals)$ forms a
chain complex, \ie, $\bdmap_{i-1}\bdmap_i = 0$.  This simple
observation is the key step to several results that follow.  To start
with, the usual homology groups $\tilde{H}_i(\Phi;\reals) = \ker
\bdmap_i / \im \bdmap_{i+1}$ are well-defined.  Recall $\betat_i(\Phi)
= \dim_{\reals}\tilde{H}_i(\Phi;\reals)$.

\begin{lemma}\label{th:betat.par}
If $\Phi$ is an interval with vertex $e$, then
$$
\betat_i(\Phi\re e) = \betat_i(\Phi)
$$
for all $i$.
\end{lemma}
\begin{proof}
First note that $\st_{\Phi}e = e * ((\Phi-e) \cap (\Phi/e)) = e *
(\Gamma, \Gamma') = (e*\Gamma, e*\Gamma')$ for some simplicial
complexes $\Gamma$ and $\Gamma'$, and so is acyclic.  Now, $\Phi$,
$\Phi \re e$, and $\st_{\Phi}e$ are all intervals, and thus chain
complexes; furthermore, by definition of $\Phi \re e$,
$$
0 \rightarrow \st_{\Phi}e \rightarrow \Phi \rightarrow \Phi \re e \rightarrow 0
$$
is a short exact sequence of chain complexes.  The resulting long
exact sequence in reduced homology (\eg, \cite[Section 24]{Munkres}), 
$$
\cdots \rightarrow \Ht_i(\st_{\Phi}e) \rightarrow \Ht_i(\Phi) \rightarrow \Ht_i(\Phi \re e) \rightarrow \Ht_{i-1}(\st_{\Phi}e) \rightarrow \cdots,
$$
becomes
$$
\cdots \rightarrow 0 \rightarrow \Ht_i(\Phi) \rightarrow \Ht_i(\Phi \re e) \rightarrow 0 \rightarrow \cdots,
$$
and the result follows immediately.
\end{proof}

We collect here the easy facts we need about how interval direct sums
and joins (and thus cones and open stars) interact with deletion,
contraction, stars, and reduction.  Each fact is either immediate from
the relevant definitions, or a routine calculation.
For the identities with the join, we assume $e$ is a vertex of $\Phi$.

\begin{align*}
(\Phi \oplus \Theta) - e     &= (\Phi - e) \oplus (\Theta - e)            & (\Phi * \Theta) - e     &= (\Phi - e) * \Theta \\ 
(\Phi \oplus \Theta)/e       &= (\Phi/e) \oplus (\Theta/e)		  & (\Phi * \Theta)/e       &= (\Phi/e) * \Theta)\\   
\st_{(\Phi \oplus \Theta)} e &= \st_{\Phi} e \oplus \st_{\Theta} e	  & \st_{(\Phi * \Theta)} e &= \st_{\Phi} e * \Theta\\
(\Phi \oplus \Theta) \re e   &= (\Phi \re e) \oplus (\Theta \re e)	  & (\Phi * \Theta) \re e   &= (\Phi \re e) * \Theta\\
\end{align*}

\section{Laplacians}\label{se:laplacians}
In this section, we define the Laplacian operators and the spectral
recursion, develop the tools we will need later to work with them, and
show that several operations on intervals, including duality
(Proposition \ref{th:recursion.dual}), preserve the property of
satisfying the spectral recursion.

\begin{definition}
The {\em {\rm (}$i$-dimensional\,{\rm )} Laplacian} of
$\Phi$ is the map $L_i(\Phi) \colon
C_i(\Phi;\reals) \rightarrow
C_i(\Phi;\reals)$ defined by
$$
L_i = L_i(\Phi) := \bdmap_{i+1}\bdmap_{i+1}^* + \bdmap_{i}^*\bdmap_{i}.
$$
\end{definition}
It is not hard to see that $L_i(\Phi)$ maps each face $[F]$ to a
linear combination of faces in $\Phi$ {\em adjacent} to $F$, that is,
faces in $\Phi$ of the form $F-v \disun w$ for some (not necessarily
distinct) vertices $v, w$, and such that $F-v \in \Phi$ or $F \disun w
\in \Phi$.  For details on the
coefficients of these linear combinations (in the simplicial complex
case, though the ideas are similar for intervals), see \cite[equations
(3.2)--(3.4)]{DuvalReiner}, but we will not need that level of
detail here.  For more information on Laplacians, also see, \eg,
\cite{Friedman,KookReinerStanton,Merris:survey}. 

Each of $\bdmap_{i+1}\bdmap_{i+1}^*$ and $\bdmap_{i}^*\bdmap_{i}$ is
positive semidefinite, since each is the composition of a linear map
and its adjoint.  Therefore, their sum $L_i$ is also positive
semidefinite, and so has only non-negative real eigenvalues.  (See
also \cite[Proposition 2.1]{Friedman}.)  These eigenvalues do not
depend on the arbitrary ordering of the vertices of $\Phi$, and are
thus invariants of $\Phi$; see, \eg, \cite[Remark 3.2]{DuvalReiner}.
Define $\s_i(\Phi)$ to be the multiset of eigenvalues of $L_i(\Phi)$,
and define $m_\lambda(L_i(\Phi))$ to be the multiplicity of $\lambda$
in $\s_i(\Phi)$.

The first result of combinatorial Hodge theory, which goes back to
Eckmann \cite{Eckmann}, is that 
\begin{equation}\label{eq:beta.m0}
m_0(L_i(\Phi))=\tilde{\beta}_i(\Phi).
\end{equation}
Though initially stated only for the case where $\Phi$ is a simplicial
complex, there is a simple proof that only relies upon $\Phi$ being a
chain complex, and so applies to all intervals $\Phi$; see
\cite[Proposition 2.1]{Friedman}.

A natural generating function for the Laplacian eigenvalues of an
interval $\Phi$ is
$$
S_{\Phi}(t,q) := \sum_{i \geq 0} t^i \sum_{\lambda \in \s_{i-1}(\Phi)} q^\lambda
        = \sum_{i, \lambda} m_\lambda(L_{i-1}(\Phi)) t^i q^\lambda.
$$
We call $S_{\Phi}$ the {\em spectrum polynomial} of $\Phi$.  It was
introduced (with slightly different indexing) for matroids in
\cite{KookReinerStanton}, and extended to relative simplicial pairs in
\cite{Duval:recursion}.  Although $S_{\Phi}$ is defined for any
interval $\Phi$, it is only truly a polynomial when the Laplacian
eigenvalues are not only non-negative, but integral as well.  This
will be true for the cases we are concerned with, primarily shifted
intervals \cite{Duval:recursion}, matroids \cite{KookReinerStanton},
and matroid pairs $(M-e,M/e)$ \cite{Duval:recursion}.

Let $F$ be a face in interval $\Phi$.  As usual, the {\em boundary} of
$F$ in $\Phi$ is the collection of faces $\{F -v \in \Phi\colon v \in
F\}$.  Similarly, the {\em coboundary} of $F$ in $\Phi$ is the
collection of faces $\{F \disun w \in \Phi\colon w \not\in F\}$.
It is not hard to see that $\bdmap_{(\Phi^*)}$ and $(\bdmap_{\Phi})^*$
each map $[F]$ to a linear combination of faces in the coboundary of
$F$ in $\Phi$.  In fact, \cite[Lemma 6.1]{Duval:recursion} states that
$\bdmap_{(\Phi^*)}$ and $(\bdmap_{\Phi})^*$
are isomorphic, up to an easy change of basis (multiplying some
basis elements by $-1$).  The easy corollary \cite[Corollary
6.2]{Duval:recursion} is that $L_i(\Phi)$ is, modulo that same change of
basis, isomorphic to $L_{n-i-2}(\Phi^*)$.  Therefore \cite[equation
(28)]{Duval:recursion},
$$
S_{\Phi^*}(t,q) = t^{\abs{E}}S_{\Phi}(t^{-1},q).
$$

By \cite[Corollary 4.3]{Duval:recursion}, 
$$
S_{\Phi * \Theta} = S_{\Phi} S_{\Theta};
$$  
it follows then that 
$$
S_{R \ostar \Phi} = t^{\abs{R}} S_{\Phi}.
$$
The following is the analogue for direct sums.  It is simpler than the
formula for disjoint union of simplicial complexes
\cite[Lemma 6.9]{Duval:recursion}, because even disjoint simplicial
complexes share the empty face.
\begin{lemma}\label{th:S.direct.sum}
If $\Phi$ and $\Theta$ are intervals such that $\Phi \oplus \Theta$ is
well-defined, then $\s_i(\Phi \oplus \Theta) = \s_i(\Phi) \cup
\s_i(\Theta)$, the multiset union of $\s_i(\Phi)$ and $\s_i(\Theta)$,
and $S_{\Phi \oplus \Theta} = S_{\Phi} + S_{\Theta}$.
\end{lemma}
\begin{proof}
Since no face in $\Theta$ is related to any face in $\Phi$, there are
no adjacencies between faces in $\Phi$ and faces in $\Theta$, nor do
any of the faces in $\Theta$ change any adjacencies in $\Phi$.
Similarly, no faces in $\Phi$ change any adjacencies in $\Theta$, and
we conclude $L_i(\Phi \oplus \Theta) = L_i(\Phi) \oplus L_i(\Theta)$.  Thus
$\s_i(\Phi \oplus \Theta) = \s_i(\Phi) \cup \s_i(\Theta)$, and so $S_{\Phi \oplus \Theta} = S_{\Phi} + S_{\Theta}$.
\end{proof}

Following \cite{DuvalReiner}, let the equivalence relation $\blambda
\circeq \bmu$ on multisets $\blambda$ and $\bmu$ denote that
$\blambda$ and $\bmu$ agree in the multiplicities of all of their {\em
non-zero parts}, \ie, that they coincide except for possibly their
number of zeros.  

\begin{lemma}\label{th:no.co.boundary}
If $\Phi$ and $\Theta$ are two intervals such that $\Phi = \Theta
\disun \calN$, where $\calN$ is a collections of faces with neither
boundary nor coboundary in $\Phi$, then $\s_i(\Phi) \circeq \s_i(\Theta)$.
\end{lemma}
\begin{proof}
Since $\Phi$ is an interval, the faces in $\calN$ are not related to
any other face in $\Phi$.  Thus $\Phi = \Theta \oplus \calN$.
Furthermore, since the faces in $\calN$ are not related to each other,
$L_i(\calN)$ is the zero matrix for all $i$, and so $\s_i(\calN)$
consists of all $0$'s.  Now apply Lemma \ref{th:S.direct.sum}.
\end{proof}

\begin{definition}
We will say that an interval $\Phi$ {\em satisfies the spectral
recursion with respect to $e$} if $e$ is a vertex of $\Phi$ and
\begin{equation}\label{eq:big}
S_\Phi(t,q) = qS_{\Phi-e}(t,q) + qtS_{\Phi/e}(t,q) + (1-q)S_{\Phi \re e}(t,q).
\end{equation}
We will say $\Phi$ {\em satisfies
the spectral recursion} if $\Phi$ satisfies the spectral recursion with respect
to every vertex in its vertex set.  (Note that Lemma \ref{th:loop.ok}
below means we need not be too particular about the vertex set of
$\Phi$.)
\end{definition}

When $\Phi$ is a simplicial complex, $\Phi \re e$ becomes
$(\Phi-e,\Phi/e)$, and equation \eqref{eq:big} immediately reduces to
the spectral recursion for simplicial complexes in
\cite{Duval:recursion}.

The statement and proof of the following lemma strongly resemble their
simplicial complex counterparts \cite[Theorem 2.4 and Corollary
4.8]{Duval:recursion}.  Here as there, specializations of the spectrum
polynomial reduce it to nice invariants of the interval, and reduce
the spectral recursion to basic recursions for those invariants.  We
sketch the proof in order to state what the spectrum polynomial and
spectral recursion reduce to in each case.

\begin{lemma}\label{th:recursion.special.values}
The spectral recursion holds for all intervals when $q=0$, $q=1$, $t=0$, or $t=-1$
\end{lemma}
\begin{proof}
If $q=0$, then by equation \eqref{eq:beta.m0}, $S_{\Phi}$ becomes
$\sum_i t^i \betat_{i-1}(\Phi)$, as in \cite[Theorem
2.4]{Duval:recursion}.  The spectral recursion then reduces to the
identity $\betat_i(\Phi) = \betat_i(\Phi \re e)$, which we established
in Lemma \ref{th:betat.par}.

If $q=1$, then $S_{\Phi}$ becomes $\sum_i t^i f_{i-1}(\Phi)$, as in
\cite[Theorem 2.4]{Duval:recursion}, where $f_i(\Phi) = \abs{\Phi_i}$.
The spectral recursion then reduces to the easy identity
\begin{equation}\label{eq:f.recursion}
f_i(\Phi) = f_i(\Phi - e) + f_{i-1}(\Phi/e).
\end{equation}

If $t=0$, then $S_{\Phi}$ becomes $q^{f_0(\Phi)}$ if $\emptyset \in
\Phi$ (as in \cite[Theorem 2.4]{Duval:recursion}), but becomes $0$
otherwise.  If $\emptyset \not\in \Phi$, then every term in the
spectral recursion becomes $0$; if, on the other hand, $\emptyset \in
\Phi$, then, as in \cite[Theorem 2.4]{Duval:recursion}, the spectral
recursion reduces to the trivial observation that $f_0(\Phi) =
f_0(\Phi -e)$ if $e$ is not a face of $\Phi$, but $f_0(\Phi) = 1+
f_0(\Phi -e)$ if $e$ is a face of $\Phi$.

If $t=-1$, then $S_{\Phi}$ becomes $\chi(\Phi) = \sum_i (-1)^i
f_i(\Phi) = \sum_i (-1)^i \betat_i(\Phi)$, the Euler characteristic of
$\Phi$, by \cite[Corollary 4.8]{Duval:recursion}.  The spectral
recursion now reduces to two easy identities about Euler
characteristic: that $\chi(\Phi) = \chi(\Phi \re e)$, which follows
from Lemma \ref{th:betat.par}; and that $\chi(\Phi) = \chi(\Phi - e) -
\chi(\Phi/e)$, which follows from the identity \eqref{eq:f.recursion}
above.
\end{proof}

If $\Phi$ is an interval and $e$ is a vertex of $\Phi$, define
$$
\calS_i(\Phi,e) = [t^i](S_{\Phi} - qS_{\Phi-e} - qtS_{\Phi/e} - (1-q)S_{\Phi \re e}),
$$
where $[t^i]p$ denotes the coefficient of $t^i$ in polynomial $p$.
Clearly, $\Phi$ satisfies the spectral recursion with respect to $e$
precisely when $\calS_i(\Phi,e)=0$ for all $i$.

\begin{lemma}\label{th:non.zeros}
Let $\Phi$ and $\Theta$ be intervals, each with vertex $e$, such that
$s_i(\Phi) \circeq \s_j(\Theta)$,
$s_i(\Phi-e) \circeq \s_j(\Theta-e)$,
$s_{i-1}(\Phi/e) \circeq \s_{j-1}(\Theta/e)$, and
$s_i(\Phi \re e) \circeq \s_j(\Theta \re e)$.
Then $\calS_i(\Phi,e) = \calS_j(\Theta,e)$.
\end{lemma}
\begin{proof}
Translating the $\circeq$ assumptions to generating functions, 
\begin{align*}
[t^i]S_{\Phi}             &= [t^j]S_{\Theta} + C_1\\
[t^i]S_{\Phi-e}           &= [t^j]S_{\Theta-e} + C_2\\
[t^{i-1}]S_{\Phi/e}       &= [t^{j-1}]S_{\Theta/e} + C_3\\
[t^i]S_{\Phi \re e} &= [t^j]S_{\Theta \re e} + C_4,
\end{align*}
where $C_1,C_2,C_3$, and $C_4$ are constants.  It is then easy to
compute
$$
\calS_i(\Phi,e) - \calS_j(\Theta,e) = (C_1 - C_4) + q(C_4 - C_2 - C_3).
$$

This makes $\calS_i(\Phi,e) - \calS_j(\Theta,e)$ a linear polynomial
in $q$.  But by Lemma \ref{th:recursion.special.values},
$\calS_i(\Phi,e) - \calS_j(\Theta,e)=0$ when $q=0$ and when $q=1$.
Therefore $\calS_i(\Phi,e) - \calS_j(\Theta,e)$ must be identically
$0$, as desired.
\end{proof}

The following two results are easy to verify directly; the third is not much harder.

\begin{lemma}\label{th:loop.ok}
If $\Phi$ is an interval and $e$ is a loop, then $\Phi$ satisfies the
spectral recursion with respect to $e$.
\end{lemma}

\begin{lemma}\label{th:vertex.ok}
The interval with only a single face, and the interval whose only two
faces are a single vertex and the empty face, each satisfy the
spectral recursion.
\end{lemma}

\begin{proposition}\label{th:recursion.dual}
Let $\Phi$ be an interval with vertex $e$.  If $\Phi$ satisfies the
spectral recursion with respect to $e$, then so does $\Phi^*$.
\end{proposition}
\begin{proof}
Calculate
\begin{align*}
S_{\Phi^*}(t,q) 
   &= t^n S_{\Phi}(t^{-1},q)\\
   &= t^n(q S_{\Phi-e}(t^{-1},q) + qt^{-1}S_{\Phi/e}(t^{-1},q)+(1-q)S_{\Phi \re e}(t^{-1},q))\\
   &= qS_{(\Phi-e)^*}(t,q) + qt^{-1}S_{(\Phi/e)^*}(t,q) + (1-q)S_{(\Phi \re e)^*}(t,q)\\
   &= qS_{e \ostar (\Phi^*/e)}(t,q) + qt^{-1}S_{e \ostar(\Phi^*-e)}(t,q) + (1-q)S_{\Phi^* \re e}(t,q)\\
   &= qtS_{\Phi^*/e}(t,q) + qS_{\Phi^*-e}(t,q) + (1-q)S_{\Phi^* \re e}(t,q).
\end{align*}
\end{proof}

Similar routine calculations establish the following two lemmas.

\begin{lemma}\label{th:recursion.direct.sum}
If $\Phi$ and $\Theta$ are intervals that satisfy the spectral
recursion with respect to $e$, and such that $\Phi \oplus \Theta$ is
well-defined, then $\Phi \oplus \Theta$ satisfies the spectral
recursion with respect to $e$.
\end{lemma}

\begin{lemma}\label{th:recursion.join}
If $\Phi$ is an interval that satisfies the spectral recursion with
respect to $e$, and $\Theta$ is another interval such that $\Phi *
\Theta$ is well-defined, then $\Phi * \Theta$ satisfies the spectral
recursion with respect to $e$.
\end{lemma}

\begin{corollary}\label{th:recursion.cone}
Let $\Phi$ be an interval.  If $\Phi$ satisfies the spectral recursion, 
then so do $v * \Phi$ and $R \ostar \Phi$.
\end{corollary}
\begin{proof}
Combine Lemmas \ref{th:vertex.ok} and \ref{th:recursion.join}
\end{proof}

\section{Skeleta}\label{se:skeleta}

The main goal of this section is to show that taking skeleta preserves
the property of satisfying the spectral recursion (Theorem
\ref{th:recursion.skeleta.hard}).  A key step is to show that skeleta
and reduction interact reasonably well (Corollary \ref{th:s.skel.par}).

\begin{definition}
We will say interval $\Phi$ is $(i,j)$-{\em dimensional} when $i \leq \dim F
\leq j$ for all $F \in \Phi$.  Note that it is {\em not} necessary for
there to be a face of every dimension between $i$ and $j$.
If $\Phi$ is an interval, we define the $(i,j)$-{\em skeleton} to be
$$
\Phi\skel{i,j} = \{F \in \Phi\colon i \leq \dim F \leq j\}
$$
\end{definition}

It is immediate that
\begin{align*}
 \Phi\skel{i,j}-e &= (\Phi - e)\skel{i,j},\\   
 \Phi\skel{i,j}/e &= (\Phi / e)\skel{i-1,j-1}.  
\end{align*}
The corresponding statement with reduction instead of deletion or
contraction is not true.  For instance, in Example \ref{ex:par}, $1256
\not\in (\Theta \re 3)\skel{1,3}$ (since $12356 \in \Theta$), but $1256
\in \Theta\skel{1,3} \re 3$ (since $12356$ is 4-dimensional, and so is
not in $\Theta\skel{1,3}$).
On the other hand, it will not be hard to show that at least
the non-zero eigenvalues of $\Phi\skel{i,j} \re e$ and $(\Phi \re
e)\skel{i,j}$ coincide.  We first need two easy technical lemmas.

\begin{lemma}\label{th:skel.par.pairs}
Let $\Phi$ be an interval with vertices $e$ and $v$.  If $F, F \disun
v \in \Phi\skel{i,j} \re e$ for some $i<j$, then $F, F \disun v
\in \Phi \re e$.
\end{lemma}
\begin{proof}
First note that $v \neq e$, since, otherwise, $F, F \disun v \in
\Phi\skel{i,j} \re e$ would be impossible.  Thus, either $e$ is
a vertex of both $F$ and $F \disun v$, or $e$ is a vertex of neither.

First assume $e \not\in F, F\disun v$.  Then $F \in \Phi\skel{i,j}
\re e$ implies $F \disun e \not\in \Phi\skel{i,j}$, and so $F
\disun e \not\in \Phi$ (note $\dim F < j$).  
But then $F \disun \{v,e\} \not\in \Phi$, since $\Phi$ is an
interval and $F \in \Phi$.
Now, with $F \disun e, F \disun \{v,e\} \not\in \Phi$, we conclude
$F, F\disun v \in \Phi \re e$.

Next assume $e \in F, F\disun v$.  Then $F \disun v \in \Phi\skel{i,j}
\re e$ implies $F \disun v - e \not\in \Phi\skel{i,j}$, and so $F
\disun v - e \not\in \Phi$ (note $\dim F \disun v > i$).  But then
$F-e \not\in \Phi$, since $\Phi$ is an interval and $F \disun v \in
\Phi$.  Now, with $F-e, F\disun v -e \not\in \Phi$, we conclude $F, F
\disun v \in \Phi \re e$.
\end{proof}

\begin{lemma}\label{th:skel.par}
Let $\Phi$ be an interval with vertex $e$.  Then
$$
\Phi\skel{i,j} \re e = (\Phi \re e)\skel{i,j} \disun \calN,
$$
where $\calN$ is a set of faces with neither boundary nor coboundary in
$\Phi\skel{i,j}\re e$.
\end{lemma}
\begin{proof}
First we show $(\Phi \re e)\skel{i,j} \subseteq \Phi\skel{i,j} \re e$.
Let $F \in (\Phi \re e)\skel{i,j}$, so $F \in \Phi \re e$ and $F \in
\Phi\skel{i,j}$.  If $e \not\in F$, then $F \disun e \not\in \Phi$, so
$F \disun e \not\in \Phi\skel{i,j}$, and so $F \in \Phi\skel{i,j} \re
e$.  If, on the other hand, $e \in F$, then $F - e \not\in \Phi$, so
$F - e \not\in \Phi\skel{i,j}$, and so $F \in \Phi\skel{i,j} \re e$.

Now let $G \in \Phi\skel{i,j} \re e$, $G \not\in (\Phi \re
e)\skel{i,j}$.  By Lemma \ref{th:skel.par.pairs}, for every $v \in G$,
we have $G-v \not\in \Phi\skel{i,j} \re e$, and, for every $w \not\in
G$, we have $G \disun w \not\in \Phi\skel{i,j} \re e$.
Therefore, $G$ has neither boundary nor coboundary in $\Phi\skel{i,j}
\re e$, as desired.
\end{proof}

\begin{corollary}\label{th:s.skel.par}
Let $\Phi$ be an interval with vertex $e$, and let $i<j$.  Then
$$
\s_k(\Phi\skel{i,j} \re e) \circeq s_k((\Phi \re e)\skel{i,j}),
$$
for all $k$.
\end{corollary}
\begin{proof}
Apply Lemma \ref{th:no.co.boundary} to Lemma \ref{th:skel.par}.
\end{proof}

The following two equations are from \cite[equation
(3.6)]{DuvalReiner}, where they are established for simplicial
complexes, but they are just easy consequences of $\Phi$ being a chain
complex.
\begin{gather}
\s_i(\Phi) \circeq \s_i(\Phi\skel{i-1,i}) \cup \s_i(\Phi\skel{i,i+1}),\label{eq:s.union.up.down}\\
\s_{i-1}(\Phi\skel{i-1,i}) \circeq \s_i(\Phi\skel{i-1,i}).\label{eq:s.up.down.same}
\end{gather}
As a result of this second equation, if $\Phi$ is $(i-1,i)$-dimensional, we will
let $\s(\Phi)$ refer to the $\circeq$ equivalence class of
$\s_{i-1}(\Phi) \circeq \s_i(\Phi)$.

\begin{lemma}\label{th:calS.up.down}
If $\Phi$ is an $(i-1,i)$-dimensional interval with vertex $e$,
then $\calS_i(\Phi,e) = \calS_{i-1}(\Phi,e)$.
\end{lemma}
\begin{proof}
By equation \eqref{eq:s.up.down.same}, since $\Phi$ is $(i-1,i)$-dimensional, $\s_{i-1}
\circeq \s_i$ for $\Phi$, $\Phi-e$, and $\Phi \re e$.  Similarly,
$\s_{i-2} \circeq \s_{i-1}$ for $\Phi/e$.  Now apply Lemma \ref{th:non.zeros}.
\end{proof}

\begin{lemma}\label{th:calS.split}
If $\Phi$ is an interval with vertex $e$, then 
$$
\calS_{i}(\Phi,e)  = \calS_{i}(\Phi\skel{i-1,i},e) + \calS_{i}(\Phi\skel{i,i+1}, e).
$$
\end{lemma}
\begin{proof}
Let $b$ and $t$ be two new vertices not in $\Phi$, and let 
$$
\Theta = (b \ostar \Phi\skel{i-1,i}) \oplus (t \ostar \Phi\skel{i,i+1}).
$$
It is immediate that $\Theta$ is well-defined, since $b \neq t$.
(Indeed, $b$ and $t$ are introduced precisely to make a direct sum of
out $\Phi\skel{i-1,i}$ and $\Phi\skel{i,i+1}$.)
It is easy to verify that 
\begin{align*}
\s_{i+1}(\Theta) &=\s_i(\Phi\skel{i-1,i}) \cup \s_i(\Phi\skel{i,i+1}) \circeq \s_i(\Phi),\\
\s_{i+1}(\Theta-e) &=\s_i((\Phi-e)\skel{i-1,i}) \cup \s_i((\Phi-e)\skel{i,i+1}) \circeq \s_i(\Phi-e),\\
\s_i(\Theta/e) &=\s_{i-1}((\Phi/e)\skel{i-2,i-1}) \cup \s_{i-1}((\Phi/e)\skel{i-1,i}) \circeq \s_{i-1}(\Phi/e),\\
\s_{i+1}(\Theta \re e) &= \s_i(\Phi\skel{i-1,i} \re e) \cup \s_i(\Phi\skel{i,i+1} \re e)\\
 &\circeq\s_i((\Phi \re e)\skel{i-1,i}) \cup \s_i((\Phi \re e)\skel{i,i+1}) \circeq \s_i(\Phi \re e);
\end{align*}
in each case, the last $\circeq$-equivalence is by equation
\eqref{eq:s.union.up.down}.  Then, by Lemma \ref{th:non.zeros},
$\calS_{i+1}(\Theta,e) = \calS_i(\Phi,e)$, and so now it is easy to
verify
\begin{align*}
\calS_{i}(\Phi,e)     &= \calS_{i+1}(\Theta,e) \nonumber\\
                      &= \calS_{i+1}(b \ostar \Phi\skel{i-1,i} \oplus t \ostar \Phi\skel{i,i+1}, e) \nonumber \\
                      &= \calS_{i+1}(b \ostar \Phi\skel{i-1,i},e) + \calS_{i+1}(t \ostar \Phi\skel{i,i+1}, e) \nonumber \\
                      &= \calS_{i}(\Phi\skel{i-1,i},e) + \calS_{i}(\Phi\skel{i,i+1}, e). 
\end{align*}
\end{proof}

\begin{lemma}\label{th:recursion.skeleta.easy}
Let $\Phi$ be an interval with vertex $e$. If every skeleton
$\Phi\skel{i-1,i}$ satisfies the
spectral recursion with respect to $e$ 
then so does $\Phi$.
\end{lemma}
\begin{proof}
This is an immediate corollary to Lemma \ref{th:calS.split}
\end{proof}

\begin{theorem}\label{th:recursion.skeleta.hard}
Let $\Phi$ be an interval with vertex $e$. If $\Phi$ satisfies the
spectral recursion with respect to $e$, then so does every skeleton
$\Phi\skel{i,j}$.
\end{theorem}
\begin{proof}
By Lemma \ref{th:recursion.skeleta.easy}, it suffices to prove that
every $\Phi\skel{i,i+1}$ satisfies the spectral recursion with respect to
$e$, which we now do by induction on $i$.

If $i \leq -2$, then $\Phi\skel{i,i+1}$ is either the interval whose
only face is the empty face, or the empty interval with no faces
whatsoever.  Either way, $\Phi\skel{i,i+1}$ trivially satisfies the
spectral recursion.

If $i > -2$, then, by induction, $\calS_i(\Phi\skel{i-1,i},e)=0$, and
by hypothesis, $\calS_i(\Phi,e)=0$.  Then by Lemma
\ref{th:calS.split}, $\calS_i(\Phi\skel{i,i+1},e)=0$, and so
$\Phi\skel{i,i+1}$ satisfies the spectral recursion with respect to
$e$, by Lemma \ref{th:calS.up.down}.
\end{proof}

\section{Shifted Intervals}\label{se:shifted}

Our main goal of this section is to show that relative simplicial
pairs that are shifted (on the same vertex order) satisfy the spectral
recursion (Theorem \ref{th:shifted.big}).  The key step is the
construction of another interval $\Phi^-$ that satisfies the spectral
recursion when $\Phi$ does; this resembles, but is more involved than,
a construction in the proof of the simplicial complex case \cite[Lemma
4.22]{Duval:recursion}.  We first translate shifted relative
simplicial pairs to shifted intervals, and show that the dual of a
shifted interval is again a shifted interval (Proposition
\ref{th:dual.shifted}).

\begin{definition}
If $F=\{f_1<\cdots<f_k\}$ and $G=\{g_1<\cdots<g_k\}$ are $k$-subsets
of integers, then $F \leq_C G$ under the {\em componentwise partial
order} if $f_p \leq g_p$ for all $p$.  A simplicial complex $\Delta$
on a vertex set of integers is {\em shifted} if
$G \leq_C H$ and $H \in \Delta$ together imply $G \in \Delta$.  An
interval $\Phi$ is {\em shifted} when $\Phi=(\Delta,\Delta')$, for
some shifted simplicial complexes $\Delta$ and $\Delta'$.
\end{definition}

We would like to replace this definition of shifted interval, which
depends on the simplicial complexes involved, to one that depends only
on the interval itself.  In order to do this, we will need a single
partial order that combines the (separate) conditions of $\Delta$
being shifted, and $\Delta$ being a simplicial complex, an idea
implicit in the work of Klivans (see \eg, \cite[Figure
1]{Klivans:shifted} or \cite[Figure 1]{Klivans:matroid}).
If $F=\{f_1<\cdots<f_k\}$ and $G=\{g_1<\cdots<g_m\}$, then $F \leq_S
G$ under the {\em shifted partial order} when $k \leq m$ and
$f_{p+m-k} \leq g_p$ for all $1 \leq p \leq k$.  In particular, it is
easy to see that if $F \subseteq G$ or $F \leq_C G$, then $F \leq_S
G$.

\begin{lemma}\label{th:shifted.simplicial.complex}
If $\Delta \subseteq 2^E$, then the following are equivalent:
\begin{enumerate}
\item\label{it:old.shifted} $\Delta$ is a shifted simplicial complex; and
\item\label{it:new.shifted} $F \leq_S H$ and $H \in \Delta$ together imply $F \in \Delta$.
\end{enumerate}
\end{lemma}
\begin{proof}
That \itref{it:new.shifted} implies \itref{it:old.shifted} is an easy
exercise.  To see that \itref{it:old.shifted} implies
\itref{it:new.shifted}, assume $\Delta$ is a shifted simplicial complex
and that $F \leq_S H \in \Delta$, and let $G$ consist of the last
$\abs{F}$ elements of $H$.  Then it is easy to see that $F \leq_C
G \subseteq H$.  Therefore $G \in \Delta$, and, consequently, $F
\in \Delta$.
\end{proof}

\begin{lemma}\label{th:alt.shift.defn}
If $\Phi \subseteq 2^E$, then the following are equivalent:
\begin{enumerate}
\item $\Phi$ is a shifted interval; and
\item $F \leq_S G \leq_S H$ and $F,H \in \Phi$ together imply $G \in \Phi$.
\end{enumerate}
\end{lemma}
\begin{proof}
Thanks to Lemma \ref{th:shifted.simplicial.complex}, the proof is
entirely analogous to that of Lemma \ref{th:interval.rel.pair}, but
with $\leq_S$ instead of $\subseteq$.
\end{proof}

The following lemma, whose easy proof is omitted, means that the
partial order $\leq_S$ is {\em admissible}.
\begin{lemma}\label{th:admissible.one}
If $v \not\in F,G$, then $F \leq_S G$ iff $F \disun v \leq_S G \disun v$.
\end{lemma}

\begin{corollary}\label{th:admissible.many}
If $A \cap F = A \cap G = \emptyset$, then $F \leq_S G$ iff $F \disun A \leq_S G \disun A$.
\end{corollary}

\begin{lemma}\label{th:dual.leqS}
If $F,G \subseteq E$, then $F \leq_S G$ iff $E-G \leq_S E-F$.
\end{lemma}
\begin{proof}
Let $A = F \cap G$ and $B=(E-F) \cap (E-G)$, and let $F'=F-A$ and
$G'=G-A$.  Thus $F=F' \disun A$, $G=G' \disun A$, $E-F=G' \disun B$,
and $E-G=F' \disun B$.  Then by Corollary \ref{th:admissible.many}
twice, $F \leq_S G$ iff $F' \leq_S G'$ iff $E-G \leq_S E-F$.
\end{proof}

\begin{proposition}\label{th:dual.shifted}
If $\Phi$ is a shifted interval, then so is $\Phi^*$.
\end{proposition}
\begin{proof}
Assume $F \leq_S G \leq_S H$ and $F,H \in \Phi^*$.  Then $E-H \leq_S
E-G \leq_S E-F$, by Lemma \ref{th:dual.leqS}, and $E-F,E-H \in \Phi$.
Therefore $E-G \in \Phi$, and so $G \in \Phi$.
\end{proof}

We have one final lemma about $\leq_S$ whose easy proof is omitted.

\begin{lemma}\label{th:mixed.leq}
If $F \leq _S G$ and $\dim F < \dim G$, then $F \disun 1 \leq_S G$.
\end{lemma}

We now turn our attention to proving that shifted intervals satisfy
the spectral recursion.  We start with a definition that does not rely
upon $\Phi$ being shifted, but which will be very useful when $\Phi$
is shifted.
If $\Phi$ is an $(i-1,i)$-dimensional interval with vertex $1$, then
define
$$
\Phi^- = \Phi - \calN_{\Phi},
$$
where 
$$
\calN_{\Phi} =  \{F \in \Phi_i\colon 1 \in F, F-1 \not\in \Phi\}
               \disun \{F \in \Phi_{i-1}\colon 1 \not\in F, F \disun 1 \not\in \Phi\}.
$$
Computing $\Phi^-$ dimension by dimension, we see that, equivalently,
\begin{align}
\Phi^- &= \{F \in \Phi_i\colon 1 \not\in F\} \disun \{F \in \Phi_i\colon 1 \in F, F-1 \in \Phi\}\nonumber\\
       &\quad \disun \{F \in \Phi_{i-1}\colon 1 \not\in F, F \disun 1 \in \Phi\} \disun \{F \in \Phi_{i-1}\colon 1 \in F\} \nonumber \\
       &= (\Phi_i - 1) \disun (1 \ostar ((\Phi_i/1) \cap (\Phi_{i-1}-1))) \nonumber\\
       &\quad \disun ((\Phi_i/1) \cap (\Phi_{i-1} -1)) \disun (1 \ostar (\Phi_{i-1}/1)) \nonumber\\
       &= (\Phi_i - 1) \disun (1*((\Phi_i/1) \cap (\Phi_{i-1} -1))) \disun (1 \ostar (\Phi_{i-1}/1)).  \label{eq:phi.minus.alt}
\end{align}

\begin{lemma}\label{th:calN.Phi}
If $\Phi$ is a shifted $(i-1,i)$-dimensional interval on vertex set
$\{1,\ldots,n\}$, then the faces of $\calN_{\Phi}$ have neither
boundary nor coboundary in $\Phi$.
\end{lemma}
\begin{proof}
Let $F \in \calN_{\Phi}$.  We split the proof into two
cases, depending on the dimension of $F$.

First assume $\dim F = i-1$.  
Then $F \in \Phi$ and $F \disun 1 \not\in \Phi$, which imply $F \disun v
\not\in \Phi$ for any $v$, since $F \subseteq F \disun 1 \leq_C F
\disun v$.
Thus, $F$ has no coboundary in $\Phi$; $F$ has no
boundary in $\Phi$ simply becasue it has minimal dimension in $\Phi$.

Now assume, on the other hand, $\dim F = i$.  
Then $F \in \Phi$ and $F-1 \not\in \Phi$, which imply $F - v \not\in
\Phi$ for any $v$, since $F - v \leq_C F - 1 \subseteq F$.
Thus,
$F$ has no boundary in $\Phi$; $F$ has no coboundary in $\Phi$ simply
because it has maximal dimension in $\Phi$.
\end{proof}

\begin{lemma}\label{th:recursion.Phi.minus}
Let $\Phi$ be a shifted $(i-1,i)$-dimensional interval on vertices
$\{1,\ldots,n\}$, and let $1 \leq e \leq n$.  Then $\Phi$ satisfies
the spectral recursion with respect to $e$ iff $\Phi^-$ does.
\end{lemma}
\begin{proof}
By Lemma \ref{th:non.zeros}, it suffices to show 
$\s(\Phi^-) \circeq \s(\Phi)$, 
$\s(\Phi^- - e) \circeq \s(\Phi - e)$, 
$\s(\Phi^- / e) \circeq \s(\Phi / e)$, and
$\s(\Phi^- \re e) \circeq \s(\Phi \re e)$.
The main tools are Lemmas \ref{th:no.co.boundary} and
\ref{th:calN.Phi}, which immediately show $\s(\Phi^-) \circeq
\s(\Phi)$.

In order to show $\s(\Phi^- \re e) \circeq \s(\Phi \re e)$, we first
claim that $\st_{\Phi^-} e = \st_{\Phi} e$.  Indeed,
$\st_{(\Phi-\calN)} e = \st_{\Phi} e$ for any set $\calN$ of faces in
$\Phi$ with neither boundary nor coboundary in $\Phi$.  Then
\begin{align*}
\s(\Phi^- \re e) 
  &= \s((\Phi-\calN_{\Phi}) - \st_{(\Phi-\calN_{\Phi})}e) 
   = \s((\Phi - \st_{\Phi}e) - \calN_{\Phi}) 
   = \s((\Phi \re e) - \calN_{\Phi})\\ 
  &\circeq \s(\Phi \re e),
\end{align*}
by Lemmas \ref{th:no.co.boundary} and \ref{th:calN.Phi}, since the
faces of $\calN_{\Phi}$ have neither boundary nor coboundary in
$\Phi$, nor in any subset of $\Phi$, such as $\Phi \re e$.

To show $\s(\Phi^- - e) \circeq \s(\Phi - e)$ and $\s(\Phi^- / e)
\circeq \s(\Phi / e)$, we split into two cases: $e=1$; and $e \neq 1$.
If $e \neq 1$, then equation \eqref{eq:phi.minus.alt} makes it easy to show that 
$\Phi^- - e = (\Phi - e)^-$ and $\Phi^- / e = (\Phi / e)^-$.
Then Lemmas \ref{th:no.co.boundary} and \ref{th:calN.Phi} show
$\s(\Phi^- - e) \circeq \s((\Phi-e)^-) \circeq \s(\Phi-e)$ and
$\s(\Phi^- / e) \circeq \s((\Phi/e)^-) \circeq \s(\Phi/e)$.

To address the $e=1$ case, first note that $\Phi^- - 1 = (\Phi -
\calN_{\Phi}) - 1 = (\Phi - 1) - (\calN_{\Phi} \cap (\Phi - 1))$.  Let
$\calN' = \calN_{\Phi} \cap (\Phi - 1)$.  Since $\calN' \subseteq
\calN_{\Phi}$, every face in $\calN'$ has neither boundary nor
coboundary in $\Phi$, nor in any subset of $\Phi$, such as $\Phi-1$.
Now apply Lemma \ref{th:no.co.boundary} to see $\s(\Phi^- - 1) \circeq
\s(\Phi - 1)$.
The proof that $\s(\Phi^-/1) \circeq \s(\Phi/1)$ proceeds similarly.
\end{proof}

\begin{definition}
Let $\Phi$ be an $(i-1,i)$-dimensional interval with vertex $1$.  Define
\begin{align*}
\Phi^+ &= \Phi^- \disun \{F \disun 1\colon 1 \not\in F, F \in \Phi_i\}
                 \disun \{F - 1\colon 1 \in F, F \in \Phi_{i-1}\}.
\end{align*}
\end{definition}

\begin{lemma}\label{th:Phi.star}
If $\Phi$ is a shifted $(i-1,i)$-dimensional interval on vertices
$\{1,\ldots,n\}$, then $\Phi^+=1*\Phi'$ for some shifted interval
$\Phi'$ on vertex set $\{2,\ldots,n\}$.
\end{lemma}
\begin{proof}
First, by equation \eqref{eq:phi.minus.alt},  
\begin{align}
\Phi^+ &= (\Phi_i - 1) \disun (1*((\Phi_i/1) \cap (\Phi_{i-1} -1))) \disun (1 \ostar (\Phi_{i-1}/1)) \nonumber\\
       &\quad   \disun (1 \ostar (\Phi_i - 1)) \disun (\Phi_{i-1}/1) \nonumber \\
       &= 1 * ((\Phi_i - 1) \disun ((\Phi_i/1) \cap (\Phi_{i-1} -1)) \disun (\Phi_{i-1}/1)).\label{eq:phi.plus.alt}
\end{align}
Now, coning preserves shiftedness of intervals, since
$1*(\Delta,\Delta') = (1*\Delta,1*\Delta')$ and, as is well-known and
easy to prove, coning preserves shiftedness of simplicial
complexes.  Equation \eqref{eq:phi.plus.alt} thus reduces the proof of this
lemma to showing that
\begin{equation}\label{eq:phi.prime}
\Phi'= (\Phi_i - 1) \disun ((\Phi_i/1) \cap (\Phi_{i-1} -1)) \disun (\Phi_{i-1}/1) 
\end{equation}
is a shifted interval. 

Equation \eqref{eq:phi.prime} means $\Phi'_i= (\Phi_i - 1)$,
$\Phi'_{i-1}= \disun ((\Phi_i/1) \cap (\Phi_{i-1} -1))$, and
$\Phi'_{i-2}= (\Phi_{i-1}/1)$, and so $G \in \Phi'$ precisely
when the following conditions are met:
\begin{enumerate}
\item\label{it:dim} $i-2 \leq \dim G \leq i$;
\item\label{it:low} if $\dim G \leq i-1$, then $G \disun 1 \in \Phi$; and
\item\label{it:high} if $\dim G \geq i-1$, then $G \in \Phi$.
\end{enumerate}
We will use the characterization of shifted intervals given in Lemma
\ref{th:alt.shift.defn} to show that $\Phi'$ is a shifted interval.
So assume $G \subseteq \{2,\ldots,n\}$; $F,H \in \Phi'$; and $F \leq_S
G \leq_S H$.  We need to show $G \in \Phi'$. 
Condition \itref{it:dim} follows directly from the hypotheses on $G$.

Next we establish condition \itref{it:low}; so assume $\dim G \leq i-1$.
First note that $\dim F \leq \dim G \leq i-1$, so $F\disun 1 \in
\Phi$, and $F \disun 1 \leq_S G \disun 1$, by Lemma \ref{th:admissible.one}.
Now, if $\dim H = \dim G \leq i-1$, then $H \disun 1 \in \Phi$ and $G
\disun 1 \leq _S H \disun 1$, by Lemma \ref{th:admissible.one}, but if $\dim H >
\dim G$, then $\dim H \geq i-1$, and so $H \in \Phi$ and, by Lemma
\ref{th:mixed.leq}, $G \disun 1 \leq_S H$.  Either way, for some
$\widetilde{H}$ (either $H$ or $H \disun 1$), $F \disun 1 \leq_S G
\disun 1 \leq_S \widetilde{H}$ and $F \disun 1, \widetilde{H} \in
\Phi$.  Thus $G \disun 1 \in \Phi$, as desired.

The proof that $G$ satisfies condition \itref{it:high} is similar; we
start by assuming $\dim G \geq i-1$.  First note that $\dim H \geq
\dim G \geq i-1$, so $H \in \Phi$ while $G \leq_S H$.  Now if $\dim F
= \dim G \geq i-1$, then $F \in \Phi$ while $F \leq_S G$, but if $\dim
F < \dim G$, then $\dim F \leq i$, and so $F \disun 1 \in \Phi$ and,
by Lemma \ref{th:mixed.leq}, $F \disun 1 \leq_S G$.  Either way, for
some $\widetilde{F}$ (either $F$ or $F \disun 1$), $\widetilde{F}
\leq_S G \leq_S H$ and $\widetilde{F}, H \in \Phi$.  Thus $G \in
\Phi$, as desired.
\end{proof}

\begin{lemma}\label{th:shifted.interval}
If $\Phi$ is a shifted $(i-1,i)$-dimensional interval, then $\Phi$
satisfies the spectral recursion.
\end{lemma}
\begin{proof}
By induction on the number of non-loop vertices.  If $\Phi$ has no
non-loop vertices, the result is trivially true.
So assume $\Phi$ has vertex set $\{1,\ldots,n\}$ with $n \geq 1$.

By Lemma \ref{th:recursion.Phi.minus}, it suffices to show $\Phi^-$
satisfies the spectral recursion.  Note that, by Lemma
\ref{th:Phi.star}, $\Phi^- = (\Phi^+)\skel{i-1,i} = (1 *
\Phi')\skel{i-1,i}$ and that $\Phi'$ is a shifted
$(i-1,i)$-dimensional interval with one less non-loop vertex (namely,
vertex 1) than $\Phi^-$, and hence fewer non-loop vertices than
$\Phi$.  By induction, then, $\Phi'$ satisfies the spectral recursion.
But since taking skeleta (Theorem \ref{th:recursion.skeleta.hard}) and
coning (Corollary \ref{th:recursion.cone}) preserve the property of
satisfying the spectral recursion, $\Phi^-$ also satisfies the
spectral recursion.
\end{proof}

\begin{theorem}\label{th:shifted.big}
If $\Phi$ is a shifted interval, then $\Phi$ satisfies the spectral recursion.
\end{theorem}
\begin{proof}
It is immediate that, since $\Phi$ is shifted, so is
$\Phi\skel{i-1,i}$ for all $i$.  By Lemma \ref{th:shifted.interval},
each $\Phi\skel{i-1,i}$ satisfies the spectral recursion.  By Lemma
\ref{th:recursion.skeleta.easy}, then, $\Phi$ satisfies the spectral
recursion.
\end{proof}

\begin{remark}
It is an easy exercise to verify that, if $\Phi$ is shifted, then so
are $\Phi-e$, $\Phi/e$, and the two direct summands of $\Phi \re e$
from Proposition \ref{th:par.direct.sum}.
\end{remark}

\section{Matroid pairs}\label{se:matroid}

In this section, we show that some matroid pairs satisfy the spectral
recursion, and conjecture that many more do as well.  We first set our
notation for matroids.  For more details, see, \eg, \cite{White}.  We
let $\calC(M)$ denote the set of circuits of matroid $M$, and $\IN(M)$
denote the independence complex, which is the simplicial complex
consisting of the independent sets of $M$, and whose Laplacian was
first studied in \cite{KookReinerStanton}.  Our notation for deletion
and contraction of intervals and simplicial complexes is consistent
with the notation for deletion and contraction of matroids, \eg,
$\IN(M-e) = \IN(M) - e$ and $\IN(M/e) = \IN(M)/e$.  Similarly, $e$ is
a loop of $M$ precisely when it is a loop of $\IN(M)$.

The existence of a strong map $N \rightarrow N'$ is the natural
condition on matroids $N$ and $N'$ to yield nice results about the
interval $(\IN(N),\IN(N'))$; see, \eg, \cite{Kung}.  Roughly speaking,
it means that the matroid structures of $N$ and $N'$ are compatible,
comparable to demanding that $\Delta$ and $\Delta'$ are shifted on the
{\em same} ordered vertex set in order for $(\Delta,\Delta')$ to be
shifted pair.  The factorization theorem (\eg, \cite[Theorem
8.2.8]{Kung}) says that one characterization of the existence of
such a strong map is that $N=M-A$ and $N'=M/A$ for some matroid $M$ with
ground set $E \disun A$.  The main result of this section is that, in
the special case where $\abs{A}=1$, \ie, $\rk N - \rk N' = \dim \IN(N)
- \dim \IN(N')=1$, the interval $(\IN(N),\IN(N'))$ satisfies the
spectral recursion.  We need first one lemma.

\begin{lemma}\label{th:matroid.decomp}
If $M$ is a matroid with ground element $e$, and $e$ is not a loop, then 
$$
(\IN(M-e),\IN(M/e)) = \bigoplus_{\condns{C \in \calC(M)}{e \in C}} (C-e) \ostar \IN(M/C).
$$
\end{lemma}
\begin{proof}
This is essentially proved in \cite[Lemmas 3.3 and
3.4]{Duval:recursion}.  We sketch the proof here, both for
completeness, and to let the language of intervals, not found in the
original, simplify some of the steps.

Let $\Phi=(\IN(M-e),\IN(M/e))$.  If $I \in \Phi$, then $I$ is
independent in $M$, but $I \disun e$ is dependent in $M$, and so there
is a {\em unique} circuit of $M$, which we denote by $\ci_M(e,I)$,
contained in $I \disun e$.  For each circuit $C \in \calC(M)$, let
$M_C = \{I \in \Phi\colon \ci_M(e,I) = C\}$.  Since each $I \in \Phi$
has a unique $\ci_M(e,I)$, the $M_C$'s partition $\Phi$.

In order to show that this partition is an interval direct sum, first
note that, if $I_1 \in M_{C_1}$ and $I_2 \in M_{C_2}$, then $I_1
\disun e$ cannot contain $C_2$, since $\ci_M(e,I)$ is the {\em unique}
circuit of $M$ contained in $I \disun e$.  Then, since $C_2 \subseteq
I_2 \disun e$, it follows that $I_2 \not\subseteq I_1$; similarly $I_1
\not\subseteq I_2$.  We conclude that all the $M_C$'s are totally
unrelated, as desired.

Finally, as in \cite{Duval:recursion}, 
\begin{align*}
M_C &= \{I \in \IN(M-e)\colon C-e \subseteq I\} \\
    &= (C-e) \ostar \IN((M-e)/(C-e)) \\
    &= (C-e) \ostar \IN(M/C).
\end{align*}
\end{proof}

\begin{theorem}\label{th:matroid.big}
If $M$ is a matroid with ground element $e$, then the matroid pair
$(\IN(M-e),\IN(M/e))$ satisfies the spectral recursion.
\end{theorem}
\begin{proof}
If $e$ is not a loop of $M$, then this is an immediate corollary to
Lemmas \ref{th:recursion.direct.sum} and \ref{th:matroid.decomp},
Corollary \ref{th:recursion.cone}, and the fact \cite[Theorem
3.18]{Duval:recursion}) that matroids satisfy the spectral recursion.
If $e$ is a loop of $M$, then it is a loop of $\IN(M)$, and so
$(\IN(M-e),\IN(M/e))=(\IN(M)-e,\IN(M)/e)=(\IN(M),\emptyset)=\IN(M)$, which satisfies the
spectral recursion.
\end{proof}

We are unable to prove anything about $(\IN(N),\IN(N'))$ if $\rk N -
\rk N' > 1$, because we don't have the analogue of Lemma
\ref{th:matroid.decomp} above.  Still, experimental evidence on
randomly chosen matroids supports the following natural conjecture.

\begin{conjecture}\label{th:all.matroids}
If there is a strong map $N \rightarrow N'$ between matroids $N$ and
$N'$, then the interval $(\IN(N),\IN(N'))$ has integral Laplacian
eigenvalues, and satisfies the spectral recursion.
\end{conjecture}

\section{Acknowledgements}
I am grateful to Vic Reiner for suggesting strong maps on matroids to me.

\newcommand{\journalname}[1]{\textit{#1}}
\newcommand{\booktitle}[1]{\textit{#1}}


\begin{thebibliography}{99}
%
\bibitem{Chari}
M. K. Chari, 
``On discrete Morse functions and combinatorial decompositions'', 
\journalname{Discrete Math.}\ {\bf 217} (2000), 101--113.
%
\bibitem{Duval:recursion}
A. M. Duval,
``A common recursion for Laplacians of
matroids and shifted simplicial complexes'',
\journalname{Doc.\ Math.},
submitted; {\tt  arXiv:math.CO/0310327}.
%
\bibitem{DuvalReiner}
A. M. Duval and V. Reiner,
``Shifted simplicial complexes are Laplacian integral'', 
\journalname{Trans.\ Amer.\ Math.\ Soc.} {\bf 354} (2002),
	 4313--4344.
%
\bibitem{Eckmann}
B. Eckmann,
``Harmonische Funktionen und Randwertaufgaben in einem Komplex'',
\journalname{Comment.\ Math.\ Helv.}\ {\bf 17} (1945), 240-255.
%
\bibitem{Forman}
R. Forman, 
``Morse theory for cell complexes'', 
\journalname{Adv.\ Math.}\ {\bf 134} (1998), 90--145.
%
\bibitem{Friedman}
J. Friedman, 
``Computing Betti numbers via combinatorial Laplacian'',
in \booktitle{Proceedings of the Twenty-eighth
Annual ACM Symposium on the Theory of Computing}, (Philadelphia, 1996), pp.\ 386--391,  
ACM, New York, 1996.
%
\bibitem{Klivans:matroid}
C. J. Klivans,
``Shifted matroid complexes'',
preprint, 2003.
%
\bibitem{Klivans:shifted}
C. J. Klivans,
``Obstructions to shiftedness'',
\journalname{Discrete Comput.\ Geom.}\ {\bf 33} (2005), 535--545.
%
\bibitem{KookReinerStanton}
W. Kook, V. Reiner, and D. Stanton,
``Combinatorial Laplacians of matroid complexes'',
\journalname{J. Amer.\ Math.\ Soc.}\ {\bf 13} (2000), 129--148.
%
\bibitem{Kung}
J. P. S. Kung,
``Strong maps'', 
in 
\cite{White},
pp.\ 224--253, 
%
\bibitem{Merris:survey}
R.~Merris,
``Laplacian matrices of graphs: a survey'',
\journalname{Lin.\ Alg.\ Appl.}\ {\bf 197/198} (1994), 143--176.
%
\bibitem{Munkres}
J. R. Munkres,
\booktitle{Elements of algebraic topology},
Addison-Wesley, Menlo Park CA, 1984.
%
\bibitem{White}
N. White (ed.),
\booktitle{Theory of matroids}, 
Encyclopedia Math.\ Appl., vol.\ 26,
Cambridge Univ.\ Press, Cambridge, 1986.

\end{thebibliography}
\end{document}